\documentclass[a4paper,12pt,dvips]{article} 
\usepackage{amssymb}
\usepackage{amscd}
\usepackage{amsmath} 
\usepackage{graphicx}
\usepackage{latexsym} 
\usepackage{enumerate}

\usepackage{theorem}
\theoremstyle{plain}
\theorembodyfont{\itshape}
\newtheorem{theorem}{Theorem}[section]

\newtheorem{lemma}[theorem]{Lemma}
 
\theorembodyfont{\rmfamily}

\def\C{\mathbb{C}}
\def\P{\mathbb{P}}
\def\H{\mathcal{H}}
\def\K{\mathcal{K}}
\def\L{\mathcal{L}} 
\def\rI{\mathrm{I}}
\def\rP{\mathrm{P}}
\def\d{\delta}  
\def\sc{\scriptstyle}
\title{\bf On an Orbifold Hamiltonian Structure \\ for the First Painlev\'e 
Equation\footnote{MSC (2010): Primary 33E17; Secondary 34M55. \,  
to appear in Journal of Mathematical Society of Japan. \,  
This work is supported by Grant-in-aid for Scientific Research, JSPS, 
No. 25400102 (C) and 23224001 (S).}} 
\author{Katsunori Iwasaki\thanks{Department of Mathematics, 
Hokkaido University, Kita 10, Nishi 8, Kita-ku, Sapporo 060-0810 Japan. 
E-mail: {\tt iwasaki@math.sci.hokudai.ac.jp}} 
and Shu Okada\thanks{Nishiwaki-Kita High School, 669-32 Gonose-cho, 
Nishiwaki 677-0014 Japan.}}
\date{} 
\begin{document}
\maketitle
\begin{abstract} 
For the first Painlev\'e equation we establish an orbifold polynomial 
Hamiltonian structure on the fibration of Okamoto's spaces and show 
that this geometric structure uniquely recovers the original Painlev\'e 
equation, thereby solving a problem posed by K. Takano.  
\\[2mm]
\end{abstract} 
\section{Introduction} \label{sec:intro} 
For each of the six Painlev\'e equations $\rP_J$, Okamoto \cite{Okamoto} 
constructed what he called the space of initial conditions. 
It is a fiber $E_t$ of a fibration $\pi : E \to T$ on which $\rP_J$ 
defines a foliation that is uniform and transversal to each fiber.  
For its construction, he first had a compact surface $\overline{E}_t$ as 
an eight-time  blowup of a Hirzebruch surface and then obtained 
$E_t = \overline{E}_t \setminus V_t$ by removing a divisor $V_t$ called 
the vertical leaves.      
Afterwards, Takano {\sl et al}. \cite{MMT,Matumiya, ST} constructed 
a symplectic atlas of $E_t$, on each chart of which $\rP_J$ enjoys 
a polynomial Hamiltonian structure, and they went on to show 
that such a structure uniquely recovers $\rP_J$.        
More precisely, they were able to do so for 
$J = \mathrm{I\!I}, \mathrm{I\!I\!I}, \mathrm{I\!V}, \mathrm{V}, 
\mathrm{V\!I}$, but left open the case $J = \mathrm{I}$.  
We settle this last case in this article.     
Independently, Chiba \cite{Chiba} solved the problem based on his 
framework of Painlev\'e equations on weighted projective spaces. 
Our approach is more classical, along the lines of Okamoto and 
Takano, where what is new for $J = \mathrm{I}$ is the consideration 
of an orbifold Hamiltonian structure. 
\par
The first Painlev\'e equation $\rP_{\mathrm{I}}$ is a nonlinear 
ordinary differential equation 
\[
\dfrac{d^2 x}{dt^2} = 6 x^2 + t,  
\]
for an unknown function $x = x(t)$ with a time variable $t \in T := \C_t$. 
If we put $y := dx/dt$ then this equation can be represented 
as a time-dependent Hamiltonian system 
\begin{equation} \label{eqn:HPI}
\dfrac{dx}{dt} = \dfrac{\partial H_{\rI}}{\partial y}, 
\quad 
\dfrac{dy}{dt} = - \dfrac{\partial H_{\rI}}{\partial x}, 
\qquad 
H_{\rI}(x,y,t) = \frac{1}{2}y^2 - 2 x^3 - t x. 
\end{equation}
\par
In order to construct the space $E_t$ for system \eqref{eqn:HPI}, 
it is sufficient to carry out an eight-time blowup of a Hirzebruch 
surface as in Okamoto \cite{Okamoto}, or alternatively a nine-time 
blowup of $\P^2$ as in Duistermaat and Joshi \cite{DJ}, followed by 
removing vertical leaves. 
But this is not sufficient for the purpose of providing a symplectic 
atlas with $E_t$.  
Indeed, for $J = \mathrm{I\!I}, \mathrm{I\!I\!I}, \mathrm{I\!V}, \mathrm{V}, 
\mathrm{V\!I}$, Takano {\sl et al}. \cite{MMT,ST} had to do some extra 
work in the course of successive blowups. 
We are in an even more intricate situation that is specific for 
$J = \mathrm{I}$.   
For this we shall carry out the following.   
\par
\begin{figure}[t]
\begin{center}
\unitlength 0.1in
\begin{picture}( 43.1600,  7.1600)(  3.5000,-13.4000)
%
\special{pn 20}%
\special{ar 408 1282 58 58  0.0000000 6.2831853}%
%
\special{pn 20}%
\special{ar 1008 1282 58 58  0.0000000 6.2831853}%
%
\special{pn 20}%
\special{ar 4608 1282 58 58  0.0000000 6.2831853}%
%
\special{pn 20}%
\special{ar 1608 1282 58 58  0.0000000 6.2831853}%
%
\special{pn 20}%
\special{ar 2208 1282 58 58  0.0000000 6.2831853}%
%
\special{pn 20}%
\special{ar 4008 1282 58 58  0.0000000 6.2831853}%
%
\special{pn 20}%
\special{ar 2808 1272 58 58  0.0000000 6.2831853}%
%
\special{pn 20}%
\special{ar 3408 1282 58 58  0.0000000 6.2831853}%
%
\special{pn 20}%
\special{pa 468 1282}%
\special{pa 958 1282}%
\special{fp}%
%
\special{pn 20}%
\special{pa 1068 1282}%
\special{pa 1558 1282}%
\special{fp}%
%
\special{pn 20}%
\special{pa 1658 1282}%
\special{pa 2148 1282}%
\special{fp}%
%
\special{pn 20}%
\special{pa 2248 1282}%
\special{pa 2738 1282}%
\special{fp}%
%
\special{pn 20}%
\special{pa 2848 1282}%
\special{pa 3338 1282}%
\special{fp}%
%
\special{pn 20}%
\special{pa 3458 1282}%
\special{pa 3948 1282}%
\special{fp}%
%
\special{pn 20}%
\special{pa 4068 1282}%
\special{pa 4558 1282}%
\special{fp}%
%
\special{pn 20}%
\special{sh 1.000}%
\special{ar 1608 682 58 58  0.0000000 6.2831853}%
%
\special{pn 20}%
\special{pa 1608 1222}%
\special{pa 1608 742}%
\special{fp}%
\end{picture}%
\end{center}
\caption{Dynkin diagram of type $E_8^{(1)}$.}
\label{fig:dynkin}
\end{figure}
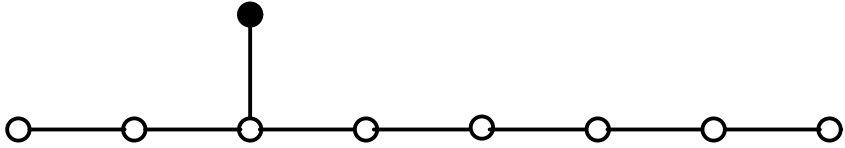

\medskip\noindent
{\bf $\bullet$ Construction of Okamoto's space}. Start with the 
Hirzebruch surface $\varSigma$ of degree $2$. 
Take a two-time blowup of $\varSigma$ to get a compact 
surface $S$, which contains a $(-2)$-curve $C$.  
Choose an open neighborhood $U$ of $C$ in $S$.     
Consider a branched double cover $(U, C) \leftarrow (V, D) 
\circlearrowleft \sigma$ ramifying along $D$, the fixed curve 
of the deck involution $\sigma$.       
Along the $(-1)$-curve $D$, take a blowdown $(V,D) \to (W,p)$ and 
let $\sigma : (W,p) \circlearrowleft$ be the induced involution.  
The result is the unique $\sigma$-fixed point $p \in W$, together 
with a pair of $\sigma$-equivalent singular points $p_{\pm} \in W$ 
of the foliation.  
To resolve the singularities $p_{\pm}$, carry out a pair of 
$\sigma$-equivariant six-time blowups $(W, p, p_{\pm}) \leftarrow 
(X, p, E_{\pm})$.  
Take a quotient $X/\sigma$, which identifies $E_+$ and $E_-$, and 
make a gluing $F = (S \setminus C) \cup (X/\sigma)$ in accordance 
with the union $S = (S \setminus C) \cup U$. 
Then $F$ is a compact surface with an $A_1$-singularity $p \in F$ 
arising from the $\sigma$-fixed point $p \in X$.  
Take a minimal resolution of $p \in F$ to obtain a smooth 
compact surface $\overline{E}_t$, which contains an $E_8^{(1)}$-type 
configuration $V_t$ of $(-2)$-curves that are the vertical leaves, 
where the black-filled node in Fig. \ref{fig:dynkin} corresponds 
to the exceptional curve for the last resolution.  
Finally we get $E_t = \overline{E}_t \setminus V_t$ 
by removing the vertical leaves $V_t$. 
Details of these processes are described in \S 2 
(see Fig. \ref{fig:blowup1} -- Fig. \ref{fig:blowup7}).     
\par\medskip\noindent
{\bf $\bullet$ Recipe for producing local charts}. We look at how a 
blowup produces two new local charts from an old one. 
Start with an $(x,y)$-plane and blow up a point $(a,0) \in \C^2_{(x,y)}$ 
on the $x$-axis $\{y = 0\}$.  
The ensuing morphism $\C^2_{(x,y)} \leftarrow \C^2_{(q,p)} \cup 
\C^2_{(Q,P)}$ is represented by @    
\[
x-a = qp, \quad y = p; \qquad x-a = Q, \quad y = QP,  
\]
where the exceptional curve is  $\{p = 0\} \cup \{Q = 0\} \cong \P^1$  
while the strict transform of $\{y = 0\}$ is $\{P = 0\}$, respectively.  
This procedure leads to a creation of two local charts 
(see Fig. \ref{fig:blowup-charts}): 
\begin{equation} \label{eqn:blowup-charts}
(x,y) \, \leadsto \, (q,p), \,\, (Q,P).    
\end{equation}
\par
Beginning with the local charts of the Hirzebruch surface $\varSigma$, we 
make a repeated application of recipe \eqref{eqn:blowup-charts} to 
produce new local charts in the course of successive blowups. 
Recall that there is one step of blowdown $(V,D) \to (W,p)$, at which 
we apply \eqref{eqn:blowup-charts} in the opposite direction. 
\begin{figure} 
\begin{center}
\unitlength 0.1in
\begin{picture}( 49.8000, 19.2000)(  2.9000,-20.6000)
%
\special{pn 20}%
\special{pa 790 1590}%
\special{pa 2390 1590}%
\special{fp}%
%
\special{pn 20}%
\special{sh 1.000}%
\special{ar 1590 1590 42 42  0.0000000 6.2831853}%
\put(14.2000,-14.7000){\makebox(0,0)[lb]{$(a,0)$}}%
\put(2.9000,-16.6000){\makebox(0,0)[lb]{$y=0$}}%
%
\special{pn 13}%
\special{pa 3380 1600}%
\special{pa 2580 1600}%
\special{fp}%
\special{sh 1}%
\special{pa 2580 1600}%
\special{pa 2648 1620}%
\special{pa 2634 1600}%
\special{pa 2648 1580}%
\special{pa 2580 1600}%
\special{fp}%
\put(27.0000,-15.0000){\makebox(0,0)[lb]{blowup}}%
%
\special{pn 20}%
\special{pa 3580 1592}%
\special{pa 5180 1592}%
\special{fp}%
%
\special{pn 20}%
\special{sh 1.000}%
\special{ar 4380 1592 42 42  0.0000000 6.2831853}%
%
\special{pn 20}%
\special{pa 4380 390}%
\special{pa 4380 1980}%
\special{fp}%
%
\special{pn 8}%
\special{ar 4380 1590 244 244  0.8960554 0.9600704}%
%
\special{pn 8}%
\special{ar 4390 1590 220 220  0.0000000 6.2831853}%
%
\special{pn 20}%
\special{sh 1.000}%
\special{ar 4380 780 42 42  0.0000000 6.2831853}%
%
\special{pn 8}%
\special{ar 4380 780 220 220  0.0000000 6.2831853}%
\put(46.3000,-19.1000){\makebox(0,0)[lb]{$(Q,P)$}}%
\put(46.8000,-8.6000){\makebox(0,0)[lb]{$(q,p)$}}%
\put(9.7000,-9.9000){\makebox(0,0)[lb]{$(x,y)$}}%
\put(52.7000,-16.5000){\makebox(0,0)[lb]{$P=0$}}%
\put(42.3000,-22.3000){\makebox(0,0)[lb]{$Q=0$}}%
\put(42.4000,-3.1000){\makebox(0,0)[lb]{$p=0$}}%
\put(26.2000,-7.9000){\makebox(0,0)[lb]{exceptional curve}}%
%
\special{pn 8}%
\special{pa 3580 850}%
\special{pa 4320 1180}%
\special{dt 0.045}%
\special{sh 1}%
\special{pa 4320 1180}%
\special{pa 4268 1136}%
\special{pa 4272 1158}%
\special{pa 4252 1172}%
\special{pa 4320 1180}%
\special{fp}%
\put(44.5000,-12.9000){\makebox(0,0)[lb]{$(-1)$}}%
\end{picture}%
\end{center}
\caption{A blowup produces two new charts from an old one.} 
\label{fig:blowup-charts}
\end{figure}
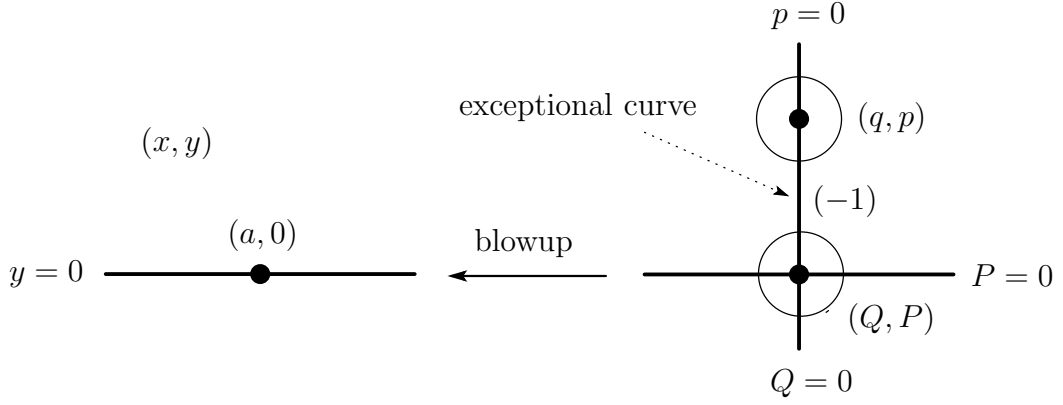
\begin{theorem} \label{thm1}
The construction mentioned above leads to the following description of 
$E_t$:    
\begin{equation} \label{eqn:Et}
E_t = \C^2_{(x,y)} \cup (\C^2_{(z,w)} \cup \C^2_{(u,v)})/\sigma,   
\end{equation}
where $\C^2_{(z,w)}$ and $\C^2_{(u,v)}$ are glued together along the 
subset $\{w \neq 0\} = \{v \neq 0 \}$ via  
\begin{equation} \label{eqn:(z,w)-(u,v)}
u = z- 2t w^{-2} - 8 w^{-6}, \qquad v = w,    
\end{equation}
with $\sigma : \C^2_{(z,w)} \cup \C^2_{(u,v)} \circlearrowleft$ being 
a holomorphic involution that restricts to  
\begin{alignat}{2} 
\sigma &: \C^2_{(z,w)} \to \C^2_{(u,v)}, \quad & (z,w) &\mapsto (u,v) = (-z,-w), \nonumber \\
\sigma &: \C^2_{(z,w)} \to \C^2_{(z,w)}, \quad & (z,w) &\mapsto (-z+2 t w^{-2}+8 w^{-6}, \, -w),  
\label{eqn:sigma2} \\
\sigma &: \C^2_{(u,v)} \to \C^2_{(u,v)}, \quad & (u,v) &\mapsto (-u-2 t v^{-2}-8 v^{-6}, \, -v),  
\nonumber
\end{alignat}
while $\C^2_{(x,y)}$ and the quotient space  
$(\C^2_{(z,w)} \cup \C^2_{(u,v)})/\sigma$ are glued together via 
\begin{alignat}{2} 
x &= \dfrac{1}{w^2}, \qquad & y &= - \dfrac{2}{w^3} - \dfrac{t w}{2} - \dfrac{w^2}{2} 
+ \dfrac{z w^3}{2},  \label{eqn:x-w} \\
x &= \dfrac{1}{v^2}, \qquad & y &= \phantom{-} \dfrac{2}{v^3} + \dfrac{t v}{2} - 
\dfrac{v^2}{2} + \dfrac{u v^3}{2},  \label{eqn:x-v}
\end{alignat}
along the subset $\{x \neq 0\} = \{w \neq 0\}/\sigma = \{v \neq 0\}/\sigma$. 
\end{theorem} 
\par
We remark that formulas \eqref{eqn:x-w} and \eqref{eqn:x-v} were already 
known to Painlev\'e \cite{Painleve} in a different context, that is, 
through the Laurent expansion around a pole of a solution, where any pole 
must be of order two so that it is converted into a simple zero via the 
transformatins $x = w^{-2} = v^{-2}$.  
\par
The total space $E$ of the fibration $\pi : E \to T$ is made up of three 
(orbifold) charts $\C^3_{(x,y,t)}$, $\C^3_{(z,w,t)}$ and $\C^3_{(u,v,t)}$ 
patched together through the symplectic mappings \eqref{eqn:(z,w)-(u,v)}, 
\eqref{eqn:x-w} and \eqref{eqn:x-v}, where by {\sl symplectic} we mean 
$\d x \wedge \d y = \d z \wedge \d w = \d u \wedge \d v$ with $\d$ being 
the relative exterior differentiation on the fibration $\pi : E \to T$ 
so that $t$ is thought of as a constant. 
In this situation one can speak of a time-dependent Hamiltonian structure 
on the fibration, which can be represented by a triple of  
Hamiltonians $H = H(x,y,t)$, $K = K(z,w,t)$ and $L = L(u,v,t)$ that 
should share a fundamental $2$-form $\Omega$ in common, to the effect that    
\[
\Omega = dy \wedge dx - dH \wedge dt = dw \wedge dz - dK \wedge dt = dv \wedge du - dL \wedge dt, 
\] 
where $d$ is the exterior differentiation on the total space $E$ so 
that $t$ is regarded as a variable. 
Under transformation rules \eqref{eqn:(z,w)-(u,v)}, \eqref{eqn:x-w} 
and \eqref{eqn:x-v}, this last condition can be written  
\begin{equation} \label{eqn:symplectic}
H = K + 1/w = L - 1/v, \qquad K = L - 2/v. 
\end{equation}
In order to speak of {\sl an orbifold Hamiltonian structure} we should 
also take into account the $\sigma$-invariance of $\Omega$. 
In view of formulas \eqref{eqn:sigma2} the condition 
$\sigma^* \Omega = \Omega$ can be written 
\begin{equation} \label{eqn:s-inv}
K \circ \sigma = K + 2/w, \qquad L \circ \sigma = L - 2/v. 
\end{equation}
\par
The first Painlev\'e equation $\rP_{\mathrm{I}}$ admits an orbifold 
Hamiltonian structure.   
Indeed, its Hamiltonian triple $\{H_{\mathrm{I}}, K_{\mathrm{I}}, 
L_{\mathrm{I}}\}$ is given by formula \eqref{eqn:HPI} together with  
\begin{equation} \label{eqn:hamiltonianPI}
\begin{split} 
K_{\mathrm{I}}(z,w,t) &= \dfrac{1}{8} w^6 z^2 - \dfrac{1}{4}(4+tw^4+w^5) z + \dfrac{1}{8} w^2(t+w)^2, \\
L_{\mathrm{I}}(u,v,t) &= \dfrac{1}{8} v^6 u^2 + \dfrac{1}{4}(4+tv^4-v^5) u + \dfrac{1}{8} v^2(t-v)^2. 
\end{split} 
\end{equation}
Note that $H_{\mathrm{I}}$, $K_{\mathrm{I}}$ and $L_{\mathrm{I}}$ are 
polynomials of their respective variables. 
Suppose that  
\begin{equation} \label{eqn:entire}
\mbox{$H$, $K$, $L$ are entire holomorphic in their respective 
variables and meromorphic on $\overline{E}$}, 
\end{equation}
where $\pi : \overline{E} \to T$ is the fibration with compactified 
fibers $\overline{E}_t$ $(t \in T)$. 
The following theorem asserts that such an orbifold Hamiltonian 
structure is unique and just coming from $\rP_{\mathrm{I}}$.     
\begin{theorem} \label{thm2}
If a function triple $\{H,K,L\}$ satisfies conditions 
\eqref{eqn:symplectic}, \eqref{eqn:s-inv} and \eqref{eqn:entire}, then 
$H = H_{\mathrm{I}}$, $K = K_{\mathrm{I}}$ and $L = L_{\mathrm{I}}$ 
modulo functions of $t \in T$.   
\end{theorem}
\par
Here we remark that a Hamiltonian makes sense only up to addition 
of a function of $t \in T$. 
Theorem \ref{thm2} is an easy consequence of the following 
function-theoretic property of $E_t$ and $\overline{E}_t$.  
\begin{theorem} \label{thm3}
Any function holomorphic on $E_t$ and meromorphic on $\overline{E}_t$ 
must be constant.  
\end{theorem}
\par
Indeed, take the differences $\H = H-H_{\mathrm{I}}$, $\K = K-K_{\mathrm{I}}$ 
and $\L = L-L_{\mathrm{I}}$. 
Since both $\{H,K,L\}$ and $\{H_{\mathrm{I}}, K_{\mathrm{I}}, 
L_{\mathrm{I}}\}$ satisfy conditions \eqref{eqn:symplectic}, 
\eqref{eqn:s-inv} and \eqref{eqn:entire}, one has $\H = \K = \L$, 
$\K \circ \sigma = \K$ and $\L \circ \sigma = \L$ so that 
$\H = \K = \L$ defines a function $h$ holomorphic on $E$ and 
meromorphic on $\overline{E}$. 
Theorem \ref{thm3} then implies that $h$ is only a function of $t \in T$. 
This proves Theorem \ref{thm2}. 
Theorem \ref{thm1} and Theorem \ref{thm3} will be proved in 
\S\ref{sec:okamoto} and \S\ref{sec:hol}, respectively.   
\section{Construction of Okamoto's Space} \label{sec:okamoto}
Our construction of $E_t$ and thus a proof of Theorem \ref{thm1} consist 
of the following twelve steps. 
\par
\begin{figure}[p]
\begin{minipage}{0.44\hsize}
\vspace{1.7mm}
\begin{center}
\unitlength 0.1in
\begin{picture}( 29.3000, 36.5000)(  1.0000,-39.4000)
%
\special{pn 13}%
\special{pa 260 3610}%
\special{pa 3030 3610}%
\special{fp}%
%
\special{pn 20}%
\special{pa 2610 520}%
\special{pa 2610 3890}%
\special{fp}%
%
\special{pn 8}%
\special{ar 610 790 192 192  0.0000000 6.2831853}%
%
\special{pn 8}%
\special{ar 624 3610 192 192  0.0000000 6.2831853}%
%
\special{pn 8}%
\special{ar 2610 780 192 192  0.0000000 6.2831853}%
%
\special{pn 8}%
\special{ar 2610 3600 192 192  0.0000000 6.2831853}%
%
\special{pn 20}%
\special{sh 1.000}%
\special{ar 612 792 50 50  0.0000000 6.2831853}%
\put(4.6000,-4.6000){\makebox(0,0)[lb]{$p_4=0$}}%
\put(6.9000,-34.0000){\makebox(0,0)[lb]{$(q_3,p_3)=(x,y)$}}%
\put(19.3000,-11.6000){\makebox(0,0)[lb]{$(q_2,p_2)$}}%
\put(22.0000,-23.2000){\makebox(0,0)[lb]{$(-2)$}}%
\put(14.7000,-7.5000){\makebox(0,0)[lb]{$(0)$}}%
\put(7.2000,-23.3000){\makebox(0,0)[lb]{$(2)$}}%
\put(14.7000,-38.6000){\makebox(0,0)[lb]{$(0)$}}%
%
\special{pn 13}%
\special{pa 620 520}%
\special{pa 620 3900}%
\special{fp}%
%
\special{pn 20}%
\special{pa 260 800}%
\special{pa 3030 800}%
\special{fp}%
\put(7.3000,-11.6000){\makebox(0,0)[lb]{$(q_4,p_4)$}}%
\put(24.8000,-4.7000){\makebox(0,0)[lb]{$p_2=0$}}%
\put(4.7000,-41.1000){\makebox(0,0)[lb]{$p_3=0$}}%
\put(24.9000,-40.9000){\makebox(0,0)[lb]{$p_1=0$}}%
\put(1.0000,-13.3000){\makebox(0,0)[lb]{$q_4=0$}}%
%
\special{pn 8}%
\special{pa 310 1160}%
\special{pa 320 850}%
\special{dt 0.045}%
\special{sh 1}%
\special{pa 320 850}%
\special{pa 298 916}%
\special{pa 318 904}%
\special{pa 338 918}%
\special{pa 320 850}%
\special{fp}%
\put(27.0000,-13.3000){\makebox(0,0)[lb]{$q_2=0$}}%
%
\special{pn 8}%
\special{pa 2910 1160}%
\special{pa 2920 850}%
\special{dt 0.045}%
\special{sh 1}%
\special{pa 2920 850}%
\special{pa 2898 916}%
\special{pa 2918 904}%
\special{pa 2938 918}%
\special{pa 2920 850}%
\special{fp}%
\put(19.6000,-34.0000){\makebox(0,0)[lb]{$(q_1,p_1)$}}%
\put(8.8000,-6.1000){\makebox(0,0)[lb]{$a^{(0)}_t$}}%
%
\special{pn 8}%
\special{pa 890 570}%
\special{pa 680 720}%
\special{dt 0.045}%
\special{sh 1}%
\special{pa 680 720}%
\special{pa 746 698}%
\special{pa 724 690}%
\special{pa 724 666}%
\special{pa 680 720}%
\special{fp}%
\end{picture}%
\end{center}
\vspace{2mm} 
\caption{Start with $\varSigma$.}
\label{fig:blowup1}
\end{minipage}
\hspace{3mm}
\begin{minipage}{0.44\hsize}
\begin{center}
\unitlength 0.1in
\begin{picture}( 28.3000, 38.2000)(  2.0000,-40.1000)
%
\special{pn 13}%
\special{pa 260 3610}%
\special{pa 3030 3610}%
\special{fp}%
%
\special{pn 20}%
\special{pa 2610 410}%
\special{pa 2610 4000}%
\special{fp}%
%
\special{pn 20}%
\special{ar 1210 2000 390 390  6.2831853 6.2831853}%
\special{ar 1210 2000 390 390  0.0000000 1.5707963}%
%
\special{pn 20}%
\special{pa 1600 2010}%
\special{pa 1600 410}%
\special{fp}%
%
\special{pn 20}%
\special{pa 270 2390}%
\special{pa 1230 2390}%
\special{fp}%
%
\special{pn 20}%
\special{pa 1240 800}%
\special{pa 3020 800}%
\special{fp}%
%
\special{pn 13}%
\special{pa 630 2020}%
\special{pa 630 4010}%
\special{fp}%
%
\special{pn 8}%
\special{ar 1600 790 192 192  0.0000000 6.2831853}%
%
\special{pn 8}%
\special{ar 622 2390 192 192  0.0000000 6.2831853}%
%
\special{pn 8}%
\special{ar 624 3610 192 192  0.0000000 6.2831853}%
%
\special{pn 8}%
\special{ar 2610 780 192 192  0.0000000 6.2831853}%
%
\special{pn 8}%
\special{ar 2610 3600 192 192  0.0000000 6.2831853}%
%
\special{pn 20}%
\special{sh 1.000}%
\special{ar 1610 800 50 50  0.0000000 6.2831853}%
\put(2.0000,-16.0000){\makebox(0,0)[lb]{$Q^{(1)}=0$}}%
\put(6.4000,-8.7000){\makebox(0,0)[lb]{$q^{(1)}=0$}}%
\put(13.9000,-3.6000){\makebox(0,0)[lb]{$p^{(1)}=0$}}%
\put(4.7000,-19.7000){\makebox(0,0)[lb]{$P^{(1)}=0$}}%
\put(7.1000,-27.5000){\makebox(0,0)[lb]{$(Q^{(1)},P^{(1)})$}}%
\put(17.0000,-11.8000){\makebox(0,0)[lb]{$(q^{(1)},p^{(1)})$}}%
\put(27.1000,-23.4000){\makebox(0,0)[lb]{$(-2)$}}%
\put(19.3000,-7.6000){\makebox(0,0)[lb]{$(-1)$}}%
\put(16.1000,-24.0000){\makebox(0,0)[lb]{$(-1)$}}%
\put(3.4000,-31.3000){\makebox(0,0)[lb]{$(1)$}}%
\put(14.6000,-35.7000){\makebox(0,0)[lb]{$(0)$}}%
\put(10.1000,-5.2000){\makebox(0,0)[lb]{$a^{(1)}_t$}}%
%
\special{pn 8}%
\special{pa 1220 490}%
\special{pa 1550 726}%
\special{dt 0.045}%
\special{sh 1}%
\special{pa 1550 726}%
\special{pa 1508 670}%
\special{pa 1508 694}%
\special{pa 1484 704}%
\special{pa 1550 726}%
\special{fp}%
%
\special{pn 8}%
\special{pa 340 1660}%
\special{pa 340 2340}%
\special{dt 0.045}%
\special{sh 1}%
\special{pa 340 2340}%
\special{pa 360 2274}%
\special{pa 340 2288}%
\special{pa 320 2274}%
\special{pa 340 2340}%
\special{fp}%
\end{picture}%
\end{center}
\caption{Blowup at $a^{(0)}_t$.}
\label{fig:blowup2}
\end{minipage}
\par\vspace{5mm}\noindent
\begin{minipage}{0.44\hsize}
\vspace{1mm}
\begin{center}
\unitlength 0.1in
\begin{picture}( 28.4000, 38.2000)(  2.1000,-40.0000)
%
\special{pn 13}%
\special{pa 270 3600}%
\special{pa 3040 3600}%
\special{fp}%
%
\special{pn 20}%
\special{pa 2610 1420}%
\special{pa 2610 4000}%
\special{fp}%
%
\special{pn 20}%
\special{ar 920 1990 390 390  6.2831853 6.2831853}%
\special{ar 920 1990 390 390  0.0000000 1.5707963}%
%
\special{pn 20}%
\special{pa 1310 2000}%
\special{pa 1310 400}%
\special{fp}%
%
\special{pn 20}%
\special{pa 260 2380}%
\special{pa 940 2380}%
\special{fp}%
\special{pa 730 2390}%
\special{pa 730 2370}%
\special{fp}%
%
\special{pn 20}%
\special{pa 910 800}%
\special{pa 2610 800}%
\special{fp}%
%
\special{pn 13}%
\special{pa 640 2010}%
\special{pa 640 4000}%
\special{fp}%
%
\special{pn 8}%
\special{ar 1310 790 192 192  0.0000000 6.2831853}%
%
\special{pn 8}%
\special{ar 632 2380 192 192  0.0000000 6.2831853}%
%
\special{pn 8}%
\special{ar 634 3600 192 192  0.0000000 6.2831853}%
%
\special{pn 8}%
\special{ar 2200 800 192 192  0.0000000 6.2831853}%
%
\special{pn 8}%
\special{ar 2620 3590 192 192  0.0000000 6.2831853}%
%
\special{pn 20}%
\special{sh 1.000}%
\special{ar 1490 800 50 50  0.0000000 6.2831853}%
\put(2.1000,-15.9000){\makebox(0,0)[lb]{$Q^{(1)}=0$}}%
\put(2.6000,-8.7000){\makebox(0,0)[lb]{$Q^{(2)}=0$}}%
\put(4.8000,-19.6000){\makebox(0,0)[lb]{$P^{(1)}=0$}}%
\put(6.9000,-28.1000){\makebox(0,0)[lb]{$(Q^{(1)},P^{(1)})$}}%
\put(23.0000,-11.6000){\makebox(0,0)[lb]{$(q^{(2)},p^{(2)})$}}%
\put(26.9000,-29.1000){\makebox(0,0)[lb]{$(-2)$}}%
\put(16.4000,-7.6000){\makebox(0,0)[lb]{$(-1)$}}%
\put(18.1000,-15.3000){\makebox(0,0)[lb]{$(-2)$}}%
\put(3.5000,-31.2000){\makebox(0,0)[lb]{$(1)$}}%
\put(14.7000,-35.6000){\makebox(0,0)[lb]{$(0)$}}%
\put(16.9000,-11.7000){\makebox(0,0)[lb]{$a^{(2)}_t$}}%
%
\special{pn 8}%
\special{pa 350 1650}%
\special{pa 350 2330}%
\special{dt 0.045}%
\special{sh 1}%
\special{pa 350 2330}%
\special{pa 370 2264}%
\special{pa 350 2278}%
\special{pa 330 2264}%
\special{pa 350 2330}%
\special{fp}%
%
\special{pn 20}%
\special{pa 2200 390}%
\special{pa 2200 1610}%
\special{fp}%
%
\special{pn 20}%
\special{pa 2400 1800}%
\special{pa 3050 1800}%
\special{fp}%
%
\special{pn 20}%
\special{ar 2390 1610 190 190  1.3156139 3.1941757}%
\put(12.3000,-3.5000){\makebox(0,0)[lb]{$P^{(2)}=0$}}%
\put(12.1000,-24.4000){\makebox(0,0)[lb]{$U$}}%
%
\special{pn 8}%
\special{pa 1710 1070}%
\special{pa 1570 880}%
\special{dt 0.045}%
\special{sh 1}%
\special{pa 1570 880}%
\special{pa 1594 946}%
\special{pa 1602 924}%
\special{pa 1626 922}%
\special{pa 1570 880}%
\special{fp}%
%
\special{pn 8}%
\special{pa 1490 800}%
\special{pa 1490 2200}%
\special{fp}%
%
\special{pn 8}%
\special{pa 1120 800}%
\special{pa 1120 2020}%
\special{fp}%
%
\special{pn 8}%
\special{pa 630 2200}%
\special{pa 970 2200}%
\special{fp}%
%
\special{pn 8}%
\special{ar 930 2010 192 192  6.2306022 6.2831853}%
\special{ar 930 2010 192 192  0.0000000 1.4659194}%
%
\special{pn 8}%
\special{pa 640 2570}%
\special{pa 1180 2570}%
\special{fp}%
%
\special{pn 8}%
\special{ar 1120 2200 372 372  6.1451793 6.2831853}%
\special{ar 1120 2200 372 372  0.0000000 1.5167948}%
%
\special{pn 8}%
\special{pa 1650 2260}%
\special{pa 1330 2120}%
\special{dt 0.045}%
\special{sh 1}%
\special{pa 1330 2120}%
\special{pa 1384 2166}%
\special{pa 1380 2142}%
\special{pa 1400 2128}%
\special{pa 1330 2120}%
\special{fp}%
\put(16.8000,-23.4000){\makebox(0,0)[lb]{$C$}}%
\put(18.6000,-23.8000){\makebox(0,0)[lb]{$(-2)$}}%
\put(3.9000,-5.9000){\makebox(0,0)[lb]{$(Q^{(2)},P^{(2)})$}}%
\put(16.7000,-29.7000){\makebox(0,0)[lb]{$S$}}%
\put(21.0000,-3.6000){\makebox(0,0)[lb]{$q^{(2)}=0$}}%
\end{picture}%
\end{center}
\vspace{2mm} 
\caption{Blowup at $a^{(1)}_t$.}
\label{fig:blowup3}
\end{minipage}
\hspace{3mm}
\begin{minipage}{0.44\hsize}
\vspace{0.7mm}
\begin{center}
\unitlength 0.1in
\begin{picture}( 28.3000, 37.9000)(  2.2000,-40.0000)
%
\special{pn 13}%
\special{pa 270 3600}%
\special{pa 3040 3600}%
\special{fp}%
%
\special{pn 20}%
\special{pa 2610 1420}%
\special{pa 2610 4000}%
\special{fp}%
%
\special{pn 20}%
\special{ar 920 1990 390 390  6.2831853 6.2831853}%
\special{ar 920 1990 390 390  0.0000000 1.5707963}%
%
\special{pn 20}%
\special{pa 1310 2000}%
\special{pa 1310 400}%
\special{fp}%
%
\special{pn 20}%
\special{pa 260 2380}%
\special{pa 940 2380}%
\special{fp}%
\special{pa 730 2390}%
\special{pa 730 2370}%
\special{fp}%
%
\special{pn 20}%
\special{pa 910 800}%
\special{pa 2610 800}%
\special{fp}%
%
\special{pn 13}%
\special{pa 640 2010}%
\special{pa 640 4000}%
\special{fp}%
%
\special{pn 8}%
\special{ar 1310 790 192 192  0.0000000 6.2831853}%
%
\special{pn 8}%
\special{ar 632 2380 192 192  0.0000000 6.2831853}%
%
\special{pn 8}%
\special{ar 634 3600 192 192  0.0000000 6.2831853}%
%
\special{pn 8}%
\special{ar 2200 800 192 192  0.0000000 6.2831853}%
%
\special{pn 8}%
\special{ar 2620 3590 192 192  0.0000000 6.2831853}%
%
\special{pn 20}%
\special{sh 1.000}%
\special{ar 1490 800 50 50  0.0000000 6.2831853}%
\put(2.2000,-16.1000){\makebox(0,0)[lb]{$s=0$}}%
\put(4.7000,-8.6000){\makebox(0,0)[lb]{$S=0$}}%
\put(5.1000,-19.7000){\makebox(0,0)[lb]{$r=0$}}%
\put(6.9000,-27.7000){\makebox(0,0)[lb]{$(r,s)$}}%
\put(23.0000,-11.6000){\makebox(0,0)[lb]{$(q^{(2)},p^{(2)})$}}%
\put(26.9000,-29.1000){\makebox(0,0)[lb]{$(-2)$}}%
\put(16.4000,-7.6000){\makebox(0,0)[lb]{$(-1)$}}%
\put(18.1000,-15.3000){\makebox(0,0)[lb]{$(-2)$}}%
\put(3.5000,-31.2000){\makebox(0,0)[lb]{$(1)$}}%
\put(14.7000,-35.6000){\makebox(0,0)[lb]{$(0)$}}%
\put(16.9000,-11.7000){\makebox(0,0)[lb]{$\tilde{a}^{(2)}_t$}}%
%
\special{pn 8}%
\special{pa 350 1650}%
\special{pa 350 2330}%
\special{dt 0.045}%
\special{sh 1}%
\special{pa 350 2330}%
\special{pa 370 2264}%
\special{pa 350 2278}%
\special{pa 330 2264}%
\special{pa 350 2330}%
\special{fp}%
%
\special{pn 20}%
\special{pa 2200 390}%
\special{pa 2200 1610}%
\special{fp}%
%
\special{pn 20}%
\special{pa 2400 1800}%
\special{pa 3050 1800}%
\special{fp}%
%
\special{pn 20}%
\special{ar 2390 1610 190 190  1.3156139 3.1941757}%
\put(15.9000,-5.2000){\makebox(0,0)[lb]{$R=0$}}%
\put(12.1000,-24.4000){\makebox(0,0)[lb]{$V$}}%
%
\special{pn 8}%
\special{pa 1710 1070}%
\special{pa 1570 880}%
\special{dt 0.045}%
\special{sh 1}%
\special{pa 1570 880}%
\special{pa 1594 946}%
\special{pa 1602 924}%
\special{pa 1626 922}%
\special{pa 1570 880}%
\special{fp}%
%
\special{pn 8}%
\special{pa 1490 800}%
\special{pa 1490 2200}%
\special{fp}%
%
\special{pn 8}%
\special{pa 1120 800}%
\special{pa 1120 2020}%
\special{fp}%
%
\special{pn 8}%
\special{pa 630 2200}%
\special{pa 970 2200}%
\special{fp}%
%
\special{pn 8}%
\special{ar 930 2010 192 192  6.2306022 6.2831853}%
\special{ar 930 2010 192 192  0.0000000 1.4659194}%
%
\special{pn 8}%
\special{pa 640 2570}%
\special{pa 1180 2570}%
\special{fp}%
%
\special{pn 8}%
\special{ar 1120 2200 372 372  6.1451793 6.2831853}%
\special{ar 1120 2200 372 372  0.0000000 1.5167948}%
%
\special{pn 8}%
\special{pa 1650 2260}%
\special{pa 1330 2120}%
\special{dt 0.045}%
\special{sh 1}%
\special{pa 1330 2120}%
\special{pa 1384 2166}%
\special{pa 1380 2142}%
\special{pa 1400 2128}%
\special{pa 1330 2120}%
\special{fp}%
\put(16.8000,-23.4000){\makebox(0,0)[lb]{$D$}}%
\put(18.6000,-23.8000){\makebox(0,0)[lb]{$(-1)$}}%
\put(7.0000,-6.6000){\makebox(0,0)[lb]{$(R,S)$}}%
%
\special{pn 20}%
\special{sh 1.000}%
\special{ar 1120 800 50 50  0.0000000 6.2831853}%
%
\special{pn 8}%
\special{pa 930 1030}%
\special{pa 1050 860}%
\special{dt 0.045}%
\special{sh 1}%
\special{pa 1050 860}%
\special{pa 996 904}%
\special{pa 1020 904}%
\special{pa 1028 926}%
\special{pa 1050 860}%
\special{fp}%
\special{pa 1050 860}%
\special{pa 1050 860}%
\special{dt 0.045}%
\put(7.5000,-11.8000){\makebox(0,0)[lb]{$\tilde{b}^{(2)}_t$}}%
%
\special{pn 8}%
\special{pa 1560 470}%
\special{pa 1340 470}%
\special{dt 0.045}%
\special{sh 1}%
\special{pa 1340 470}%
\special{pa 1408 490}%
\special{pa 1394 470}%
\special{pa 1408 450}%
\special{pa 1340 470}%
\special{fp}%
\put(12.4000,-3.9000){\makebox(0,0)[lb]{$\curvearrowright$}}%
\put(11.0000,-3.8000){\makebox(0,0)[lb]{$\sigma$}}%
\end{picture}%
\end{center}
\vspace{2mm}
\caption{Double cover.}
\label{fig:blowup4}
\end{minipage}
\end{figure} 
{\bf 1.} The Hirzebruch surface $\varSigma$ of degree $2$ is made up 
of four local charts $\C^2_{(q_i,p_i)}$, $i=1,2,3,4$, glued together 
according to the relations: 
\begin{equation} \label{eqn:hirzebruch}
q_1 q_2 = 1, \quad p_1 = -q_2 ^2 p_2; \qquad q_3 = q_1, \quad p_1 p_3 = 1; 
\qquad q_4 = q_2, \quad p_2 p_4 = 1,  
\end{equation}
where $(q_3,p_3) = (x,y)$ is the original chart for system \eqref{eqn:HPI}. 
Consider the Pfaffian system on $\C^2_{(x,y)} \times T$ defined by 
formula \eqref{eqn:HPI} and extend it to the entire space $\varSigma \times T$. 
For each $t \in T$ the associated foliation has two vertical leaves  
$\{p_1 = 0\} \cup \{p_2 = 0\} \cong \P^1$ and 
$\{q_2 = 0\} \cup \{q_4 = 0\} \cong \P^1$, together with an accessible singular 
point $a^{(0)}_t = \{(q_4, p_4) = (0, 0) \}$ (see Fig. \ref{fig:blowup1}). 
In what follows by a singularity we always mean an accessible singularity.  
\par
{\bf 2.} Blowup at $a^{(0)}_t$ produces two new charts $(q^{(1)}, p^{(1)})$ 
and $(Q^{(1)}, P^{(1)})$ such that    
\[
q_4 = q^{(1)}p^{(1)}, \quad p_4=p^{(1)}; \qquad q_4=Q^{(1)}, \quad p_4 = Q ^{(1)} P^{(1)}.  
\]
Rewrite the Pfaffian system in terms of the new charts. 
The ensuing foliation has two vertical leaves; the exceptional curve   
$\{ p^{(1)} = 0\} \cup \{Q^{(1)}=0\}$ of the blowup and the proper image 
$\{q^{(1)} = 0\}$ of $\{q_4=0\}$, together with a singular point 
$a^{(1)}_t = \{(q^{(1)}, p^{(1)})=(0,0)\}$ (see Fig. \ref{fig:blowup2}). 
\par
{\bf 3.} Blowup at $a^{(1)}_t$ produces two new charts $(q^{(2)}, p^{(2)})$ 
and $(Q^{(2)}, P^{(2)})$ such that 
\[
q^{(1)} = q^{(2)}p^{(2)}, \quad p^{(1)} = p^{(2)}; \qquad 
q^{(1)} = Q^{(2)}, \quad p^{(1)} = Q^{(2)} P^{(2)}. 
\]
In terms of the new charts there are three vertical leaves; the exceptional curve   
$\{p^{(2)} = 0\} \cup \{Q^{(2)} = 0\}$ of the blowup, the proper images    
$\{q^{(2)}=0\}$ of $\{q^{(1)} = 0\}$ and $\{P^{(2)}=0\}$ of $\{p^{(1)}=0\}$, 
together with a singular point $a^{(2)}_t = \{(Q^{(2)}, P^{(2)}) = (0,4)\}$ 
(see Fig. \ref{fig:blowup3}).   
\par
{\bf 4.} Consider the $(-2)$-curve $C = \{Q^{(1)} = 0\} \cup \{P^{(2)} = 0\}$ and 
its tubular neighborhood $U = \C^2_{(Q^{(1)}, P^{(1)})} \cup \C^2_{(Q^{(2)}, P^{(2)})}$. 
Let $(U, C) \leftarrow (V, D)$ with $V = \C^2_{(r, s)} \cup \C^2_{(R, S)}$ be 
the branched double covering ramifying along $D = \{s = 0\} \cup \{R = 0\}$, 
which is defined by
\[
Q^{(1)} = s^2, \quad P^{(1)} = r; \qquad Q^{(2)} = S, \quad P^{(2)} = R^2. 
\]
The deck involution $\sigma : V \circlearrowleft$ maps    
$(r, s) \mapsto (r, -s)$ on $\C^2_{(r,s)}$ and $(R, S) \mapsto (-R, S)$ 
on $\C ^2_{(R, S)}$, respectively. 
Three vertical leaves mentioned in step 3 become $\{S=0\}$, $\{q^{(2)} = 0\}$ 
and $\{R=0\}$, respectively, while the singular point $a_t^{(2)} \in U$ lifts 
up to a pair of $\sigma$-equivalent points in $V$: 
\begin{equation} \label{eqn:ab-tilde}
\tilde{a}^{(2)}_t = \{(R, S) = (2, 0) \} 
\,\, \oplus \,\,      
\tilde{b}^{(2)}_t = \{(R, S) = (-2, 0)\},  
\end{equation}
where $(*) \oplus (**)$ indicates that $(*)$ and $(**)$ are permuted by 
the involution $\sigma$ (see Fig. \ref{fig:blowup4}). 
In what follows $(*) \oplus (**)$ will be thought of as a single 
(that is, not a dual) object.   
\par
{\bf 5.} Since the branching locus $C$ downstairs is a $(-2)$-curve, 
the ramifying locus $D$ upstairs is a $(-1)$-curve that can be blown down 
into a smooth point $p$.   
Blowdown of $D$ into $p$ induces a morphism $(V,D) \to (W,p)$ with 
$W = \C^2_{(r^{(2)},s^{(2)})}$ and $p = \{(r^{(2)},s^{(2)}) = (0,0)\}$ such that  
\[
r^{(2)} = r s, \quad s^{(2)} = s; \qquad r^{(2)} = R, \quad s^{(2)} = R S. 
\]
\begin{figure}[p]
\begin{minipage}{0.44\hsize}
\begin{center}
\unitlength 0.1in
\begin{picture}( 27.7000, 36.0000)(  2.6000,-38.1000)
%
\special{pn 13}%
\special{pa 260 3410}%
\special{pa 3030 3410}%
\special{fp}%
%
\special{pn 20}%
\special{pa 2610 210}%
\special{pa 2610 3800}%
\special{fp}%
%
\special{pn 20}%
\special{ar 1570 1800 390 390  6.2831853 6.2831853}%
\special{ar 1570 1800 390 390  0.0000000 1.5707963}%
%
\special{pn 20}%
\special{pa 1960 1810}%
\special{pa 1960 210}%
\special{fp}%
%
\special{pn 20}%
\special{pa 260 2190}%
\special{pa 1590 2190}%
\special{fp}%
%
\special{pn 20}%
\special{pa 1600 600}%
\special{pa 3020 600}%
\special{fp}%
%
\special{pn 13}%
\special{pa 990 1820}%
\special{pa 990 3810}%
\special{fp}%
%
\special{pn 8}%
\special{ar 990 3410 192 192  0.0000000 6.2831853}%
%
\special{pn 20}%
\special{sh 1.000}%
\special{ar 1292 2192 50 50  0.0000000 6.2831853}%
\put(27.1000,-21.4000){\makebox(0,0)[lb]{$(-2)$}}%
\put(20.6000,-5.5000){\makebox(0,0)[lb]{$(-2)$}}%
\put(20.2000,-16.4000){\makebox(0,0)[lb]{$(0)$}}%
\put(7.0000,-30.6000){\makebox(0,0)[lb]{$(2)$}}%
\put(15.8000,-33.6000){\makebox(0,0)[lb]{$(0)$}}%
\put(12.3000,-21.0000){\makebox(0,0)[lb]{$p_+$}}%
%
\special{pn 20}%
\special{sh 1.000}%
\special{ar 650 2180 50 50  0.0000000 6.2831853}%
\put(5.3000,-20.9000){\makebox(0,0)[lb]{$p_-$}}%
%
\special{pn 8}%
\special{ar 980 2180 194 194  2.8587446 6.2831853}%
\special{ar 980 2180 194 194  0.0000000 1.8547980}%
\put(9.2000,-18.2000){\makebox(0,0)[lb]{$\curvearrowright$}}%
\put(10.7000,-25.6000){\makebox(0,0)[lb]{$W \, (r^{(2)},s^{(2)})$}}%
\put(9.5000,-16.6000){\makebox(0,0)[lb]{$\sigma$}}%
\put(2.8000,-14.8000){\makebox(0,0)[lb]{$r^{(2)}=0$}}%
%
\special{pn 8}%
\special{pa 400 1520}%
\special{pa 930 1900}%
\special{dt 0.045}%
\special{sh 1}%
\special{pa 930 1900}%
\special{pa 888 1846}%
\special{pa 888 1870}%
\special{pa 864 1878}%
\special{pa 930 1900}%
\special{fp}%
\put(2.8000,-27.1000){\makebox(0,0)[lb]{$s^{(2)}=0$}}%
%
\special{pn 8}%
\special{pa 320 2510}%
\special{pa 320 2220}%
\special{dt 0.045}%
\special{sh 1}%
\special{pa 320 2220}%
\special{pa 300 2288}%
\special{pa 320 2274}%
\special{pa 340 2288}%
\special{pa 320 2220}%
\special{fp}%
%
\special{pn 20}%
\special{sh 1.000}%
\special{ar 990 2190 50 50  0.0000000 6.2831853}%
\put(8.1000,-23.6000){\makebox(0,0)[lb]{$p$}}%
\end{picture}%
\end{center}
\caption{Blowdown of $D$.}
\label{fig:blowup5}
\end{minipage}
\hspace{3mm}
\begin{minipage}{0.44\hsize}
\begin{center}
\unitlength 0.1in
\begin{picture}( 27.7000, 36.0000)( 14.4000,-37.9000)
%
\special{pn 13}%
\special{pa 1440 3390}%
\special{pa 4210 3390}%
\special{fp}%
%
\special{pn 20}%
\special{pa 3790 190}%
\special{pa 3790 3780}%
\special{fp}%
%
\special{pn 20}%
\special{ar 2750 1780 390 390  6.2831853 6.2831853}%
\special{ar 2750 1780 390 390  0.0000000 1.5707963}%
%
\special{pn 20}%
\special{pa 3140 1790}%
\special{pa 3140 190}%
\special{fp}%
%
\special{pn 20}%
\special{pa 1440 2170}%
\special{pa 2770 2170}%
\special{fp}%
%
\special{pn 20}%
\special{pa 2780 580}%
\special{pa 4200 580}%
\special{fp}%
%
\special{pn 13}%
\special{pa 2170 1800}%
\special{pa 2170 3790}%
\special{fp}%
%
\special{pn 8}%
\special{ar 2470 2160 192 192  0.0000000 6.2831853}%
%
\special{pn 8}%
\special{ar 2170 3390 192 192  0.0000000 6.2831853}%
%
\special{pn 8}%
\special{ar 1830 2170 192 192  0.0000000 6.2831853}%
%
\special{pn 20}%
\special{sh 1.000}%
\special{ar 2470 1510 50 50  0.0000000 6.2831853}%
\put(38.9000,-21.2000){\makebox(0,0)[lb]{$(-2)$}}%
\put(32.4000,-5.3000){\makebox(0,0)[lb]{$(-2)$}}%
\put(32.0000,-16.2000){\makebox(0,0)[lb]{$(-1)$}}%
\put(19.0000,-31.2000){\makebox(0,0)[lb]{$(2)$}}%
\put(27.6000,-33.4000){\makebox(0,0)[lb]{$(0)$}}%
\put(27.2000,-18.9000){\makebox(0,0)[lb]{$a^{(3)}_t$}}%
%
\special{pn 20}%
\special{sh 1.000}%
\special{ar 1830 1500 50 50  0.0000000 6.2831853}%
\put(14.5000,-19.3000){\makebox(0,0)[lb]{$b^{(3)}_t$}}%
%
\special{pn 20}%
\special{pa 1830 2570}%
\special{pa 1840 880}%
\special{fp}%
%
\special{pn 20}%
\special{pa 2470 2560}%
\special{pa 2470 880}%
\special{fp}%
%
\special{pn 8}%
\special{ar 2470 1510 192 192  0.0000000 6.2831853}%
%
\special{pn 8}%
\special{ar 1830 1500 192 192  0.0000000 6.2831853}%
\put(14.9000,-27.8000){\makebox(0,0)[lb]{$U^{(3)}=0$}}%
\put(23.7000,-27.8000){\makebox(0,0)[lb]{$Z^{(3)}=0$}}%
%
\special{pn 20}%
\special{sh 1.000}%
\special{ar 2170 2170 50 50  0.0000000 6.2831853}%
\put(20.4000,-23.5000){\makebox(0,0)[lb]{$p$}}%
\put(23.1000,-8.4000){\makebox(0,0)[lb]{$w^{(3)}=0$}}%
\put(15.7000,-8.5000){\makebox(0,0)[lb]{$v^{(3)}=0$}}%
%
\special{pn 8}%
\special{pa 1650 1750}%
\special{pa 1760 1560}%
\special{dt 0.045}%
\special{sh 1}%
\special{pa 1760 1560}%
\special{pa 1710 1608}%
\special{pa 1734 1606}%
\special{pa 1744 1628}%
\special{pa 1760 1560}%
\special{fp}%
\special{pa 1760 1560}%
\special{pa 1760 1560}%
\special{dt 0.045}%
%
\special{pn 8}%
\special{pa 2690 1750}%
\special{pa 2550 1560}%
\special{dt 0.045}%
\special{sh 1}%
\special{pa 2550 1560}%
\special{pa 2574 1626}%
\special{pa 2582 1604}%
\special{pa 2606 1602}%
\special{pa 2550 1560}%
\special{fp}%
\put(25.2000,-12.9000){\makebox(0,0)[lb]{$\sc{(z^{(3)},w^{(3)})}$}}%
\put(18.8000,-12.8000){\makebox(0,0)[lb]{$\sc{(u^{(3)},v^{(3)})}$}}%
\put(21.4000,-6.9000){\makebox(0,0)[lb]{$\curvearrowright$}}%
\put(21.6000,-5.8000){\makebox(0,0)[lb]{$\sigma$}}%
\end{picture}%
\end{center}
\caption{Blowup at $p_+ \oplus p_-$.}
\label{fig:blowup6}
\end{minipage}
\vspace{5mm}
\begin{center}
\unitlength 0.1in
\begin{picture}( 55.4000, 39.6000)(  1.7000,-41.2000)
%
\special{pn 13}%
\special{pa 180 3790}%
\special{pa 5700 3800}%
\special{fp}%
%
\special{pn 20}%
\special{pa 820 2990}%
\special{pa 4010 2990}%
\special{fp}%
%
\special{pn 20}%
\special{pa 178 610}%
\special{pa 1258 610}%
\special{da 0.070}%
%
\special{pn 20}%
\special{pa 980 390}%
\special{pa 980 1650}%
\special{fp}%
%
\special{pn 20}%
\special{pa 650 1410}%
\special{pa 1770 1410}%
\special{fp}%
%
\special{pn 20}%
\special{pa 1146 2200}%
\special{pa 2234 2200}%
\special{fp}%
%
\special{pn 20}%
\special{pa 1450 1120}%
\special{pa 1450 2490}%
\special{fp}%
%
\special{pn 20}%
\special{pa 1940 2000}%
\special{pa 1940 3200}%
\special{fp}%
%
\special{pn 13}%
\special{pa 2410 2700}%
\special{pa 2410 4120}%
\special{fp}%
%
\special{pn 20}%
\special{pa 4642 620}%
\special{pa 3562 620}%
\special{da 0.070}%
%
\special{pn 20}%
\special{pa 3860 380}%
\special{pa 3860 1650}%
\special{fp}%
%
\special{pn 20}%
\special{pa 4170 1420}%
\special{pa 3050 1420}%
\special{fp}%
%
\special{pn 20}%
\special{pa 3674 2210}%
\special{pa 2586 2210}%
\special{fp}%
%
\special{pn 20}%
\special{pa 3370 1110}%
\special{pa 3370 2500}%
\special{fp}%
%
\special{pn 20}%
\special{pa 4990 800}%
\special{pa 4990 2430}%
\special{fp}%
%
\special{pn 20}%
\special{ar 4506 2380 488 610  6.2831853 6.2831853}%
\special{ar 4506 2380 488 610  0.0000000 1.5707963}%
%
\special{pn 20}%
\special{pa 5470 780}%
\special{pa 5470 4100}%
\special{fp}%
%
\special{pn 20}%
\special{pa 4740 1010}%
\special{pa 5710 1010}%
\special{fp}%
%
\special{pn 20}%
\special{sh 1.000}%
\special{ar 2410 2990 50 50  0.0000000 6.2831853}%
%
\special{pn 8}%
\special{ar 980 610 136 136  0.0000000 0.0882353}%
\special{ar 980 610 136 136  0.3529412 0.4411765}%
\special{ar 980 610 136 136  0.7058824 0.7941176}%
\special{ar 980 610 136 136  1.0588235 1.1470588}%
\special{ar 980 610 136 136  1.4117647 1.5000000}%
\special{ar 980 610 136 136  1.7647059 1.8529412}%
\special{ar 980 610 136 136  2.1176471 2.2058824}%
\special{ar 980 610 136 136  2.4705882 2.5588235}%
\special{ar 980 610 136 136  2.8235294 2.9117647}%
\special{ar 980 610 136 136  3.1764706 3.2647059}%
\special{ar 980 610 136 136  3.5294118 3.6176471}%
\special{ar 980 610 136 136  3.8823529 3.9705882}%
\special{ar 980 610 136 136  4.2352941 4.3235294}%
\special{ar 980 610 136 136  4.5882353 4.6764706}%
\special{ar 980 610 136 136  4.9411765 5.0294118}%
\special{ar 980 610 136 136  5.2941176 5.3823529}%
\special{ar 980 610 136 136  5.6470588 5.7352941}%
\special{ar 980 610 136 136  6.0000000 6.0882353}%
%
\special{pn 8}%
\special{ar 2400 3800 136 136  0.0000000 6.2831853}%
%
\special{pn 8}%
\special{ar 3840 620 136 136  0.0000000 0.0882353}%
\special{ar 3840 620 136 136  0.3529412 0.4411765}%
\special{ar 3840 620 136 136  0.7058824 0.7941176}%
\special{ar 3840 620 136 136  1.0588235 1.1470588}%
\special{ar 3840 620 136 136  1.4117647 1.5000000}%
\special{ar 3840 620 136 136  1.7647059 1.8529412}%
\special{ar 3840 620 136 136  2.1176471 2.2058824}%
\special{ar 3840 620 136 136  2.4705882 2.5588235}%
\special{ar 3840 620 136 136  2.8235294 2.9117647}%
\special{ar 3840 620 136 136  3.1764706 3.2647059}%
\special{ar 3840 620 136 136  3.5294118 3.6176471}%
\special{ar 3840 620 136 136  3.8823529 3.9705882}%
\special{ar 3840 620 136 136  4.2352941 4.3235294}%
\special{ar 3840 620 136 136  4.5882353 4.6764706}%
\special{ar 3840 620 136 136  4.9411765 5.0294118}%
\special{ar 3840 620 136 136  5.2941176 5.3823529}%
\special{ar 3840 620 136 136  5.6470588 5.7352941}%
\special{ar 3840 620 136 136  6.0000000 6.0882353}%
%
\special{pn 8}%
\special{ar 440 600 136 136  0.0000000 6.2831853}%
%
\special{pn 8}%
\special{ar 4370 610 136 136  0.0000000 6.2831853}%
\put(50.5000,-30.0000){\makebox(0,0)[lb]{$(-2)$}}%
\put(45.6000,-22.6000){\makebox(0,0)[lb]{$(-1)$}}%
\put(50.5000,-9.6000){\makebox(0,0)[lb]{$(-2)$}}%
\put(10.2000,-13.8000){\makebox(0,0)[lb]{$(-2)$}}%
\put(34.3000,-13.8000){\makebox(0,0)[lb]{$(-2)$}}%
\put(34.4000,-19.2000){\makebox(0,0)[lb]{$(-2)$}}%
\put(10.5000,-18.9000){\makebox(0,0)[lb]{$(-2)$}}%
\put(15.4000,-27.6000){\makebox(0,0)[lb]{$(-2)$}}%
\put(29.5000,-27.6000){\makebox(0,0)[lb]{$(-2)$}}%
\put(15.3000,-24.2000){\makebox(0,0)[lb]{$(-2)$}}%
\put(29.7000,-24.2000){\makebox(0,0)[lb]{$(-2)$}}%
\put(39.1000,-11.1000){\makebox(0,0)[lb]{$(-2)$}}%
\put(5.8000,-11.2000){\makebox(0,0)[lb]{$(-2)$}}%
\put(22.4000,-31.9000){\makebox(0,0)[lb]{$p$}}%
\put(13.2000,-6.6000){\makebox(0,0)[lb]{$U^{(8)}=0$}}%
\put(8.9000,-3.4000){\makebox(0,0)[lb]{$V^{(8)}=0$}}%
\put(28.8000,-6.4000){\makebox(0,0)[lb]{$Z^{(8)}=0$}}%
\put(33.6000,-3.3000){\makebox(0,0)[lb]{$W^{(8)}=0$}}%
\put(1.8000,-4.2000){\makebox(0,0)[lb]{$(u,v)$}}%
\put(42.5000,-4.4000){\makebox(0,0)[lb]{$(z,w)$}}%
\put(25.0000,-40.7000){\makebox(0,0)[lb]{$(x,y)$}}%
\put(22.8000,-17.9000){\makebox(0,0)[lb]{$\longleftrightarrow$}}%
\put(23.9000,-16.6000){\makebox(0,0)[lb]{$\sigma$}}%
%
\special{pn 8}%
\special{pa 2690 3430}%
\special{pa 2490 3100}%
\special{dt 0.045}%
\special{sh 1}%
\special{pa 2490 3100}%
\special{pa 2508 3168}%
\special{pa 2518 3146}%
\special{pa 2542 3148}%
\special{pa 2490 3100}%
\special{fp}%
\put(26.3000,-35.9000){\makebox(0,0)[lb]{fixed point of $\sigma$}}%
\put(47.0000,-6.5000){\makebox(0,0)[lb]{$w=0$}}%
\put(1.7000,-12.7000){\makebox(0,0)[lb]{$v=0$}}%
%
\special{pn 8}%
\special{pa 210 1130}%
\special{pa 210 650}%
\special{dt 0.045}%
\special{sh 1}%
\special{pa 210 650}%
\special{pa 190 718}%
\special{pa 210 704}%
\special{pa 230 718}%
\special{pa 210 650}%
\special{fp}%
%
\special{pn 20}%
\special{pa 4280 2990}%
\special{pa 4530 2990}%
\special{fp}%
%
\special{pn 20}%
\special{pa 2880 1990}%
\special{pa 2880 3190}%
\special{fp}%
%
\special{pn 8}%
\special{pa 4290 2410}%
\special{pa 4060 3320}%
\special{dt 0.045}%
\special{pa 4060 3320}%
\special{pa 4060 3320}%
\special{dt 0.045}%
\special{pa 4060 3320}%
\special{pa 4060 3320}%
\special{dt 0.045}%
%
\special{pn 8}%
\special{pa 520 2390}%
\special{pa 750 3300}%
\special{dt 0.045}%
\special{pa 750 3300}%
\special{pa 750 3300}%
\special{dt 0.045}%
\special{pa 750 3300}%
\special{pa 750 3300}%
\special{dt 0.045}%
%
\special{pn 8}%
\special{pa 760 3310}%
\special{pa 4050 3320}%
\special{dt 0.045}%
\put(30.0000,-18.2000){\makebox(0,0)[lb]{$E_+$}}%
\put(16.3000,-18.1000){\makebox(0,0)[lb]{$E_-$}}%
\put(23.4000,-11.9000){\makebox(0,0)[lb]{$X$}}%
\end{picture}%
\end{center}
\caption{A ``deer" whose dual ``horns" $E_{\pm}$ are identified by $\sigma$.}
\label{fig:blowup7}
\end{figure} 
The induced involution $\sigma : (W,p) \circlearrowleft$ maps  
$(r^{(2)}, s^{(2)}) \mapsto (-r^{(2)}, -s^{(2)})$. 
Through the blowdown morphism the vertical leaf $\{S = 0\}$ descends 
to $\{s^{(2)} = 0 \}$, while the singular point \eqref{eqn:ab-tilde} to 
\[ 
p_+ = \{(r^{(2)}, s^{(2)}) = (2, 0) \} 
\,\, \oplus \,\,     
p_- = \{(r^{(2)}, s^{(2)}) = (-2, 0) \}. 
\qquad \mbox{(see Fig. \ref{fig:blowup5})} 
\]
\par
{\bf 6.} Blowup at $p_+ \oplus p_-$ produces new charts 
$(z^{(3)}, w^{(3)}) \oplus (u^{(3)}, v^{(3)})$ and 
$(Z^{(3)}, W^{(3)}) \oplus (U^{(3)}, V^{(3)})$ such that 
\begin{alignat*}{4}
r^{(2)} &= \phantom{-}2+z^{(3)}w^{(3)}, \quad & s^{(2)} &= w ^{(3)}; \qquad & 
r^{(2)} &= \phantom{-}2+Z ^{(3)}, \quad & s^{(2)} &= Z^{(3)}W^{(3)}, \\
r^{(2)} &= -2+u^{(3)}v^{(3)}, \quad & s^{(2)} &= v ^{(3)}; \qquad & 
r^{(2)} &= -2+U ^{(3)}, \quad & s^{(2)} &= U^{(3)}V^{(3)},  
\end{alignat*}
where the induced involution $\sigma$ maps $(z^{(3)}, w^{(3)}) \mapsto 
(z^{(3)} + 4/w^{(3)}, \, -w^{(3)})$ on $\C ^2 _{(z^{(3)},w^{(3)})}$ and  
$(u^{(3)}, v^{(3)}) \mapsto (u^{(3)} - 4/v^{(3)}, \, -v^{(3)})$ on 
$\C ^2 _{(u^{(3)},v^{(3)})}$, respectively.  
In terms of the new charts there are two vertical leaves; the exceptional 
curve $\{w^{(3)} = 0\} \cup \{Z^{(3)} = 0\} \oplus \{v^{(3)} = 0\} \cup \{ U^{(3)} = 0\}$ 
and the proper image $\{W^{(3)} = 0\} = \{V^{(3)} = 0\}$ of $\{s^{(2)} = 0\}$, 
together with a singular point $a^{(3)}_t = \{(z^{(3)}, w^{(3)}) = (0, 0)\} \oplus 
b^{(3)}_t = \{(u^{(3)}, v^{(3)}) = (0, 0)\}$ (see Fig. \ref{fig:blowup6}). 
\par
{\bf 7.} Blowup at $a^{(3)}_t \oplus b^{(3)}_t$ produces new charts 
$(z^{(4)}, w^{(4)}) \oplus (u^{(4)}, v^{(4)})$  and 
$(Z^{(4)}, W^{(4)}) \oplus (U^{(4)}, V^{(4)})$ such that 
\begin{alignat*}{4}
z^{(3)} &= z^{(4)}w^{(4)}, \quad & w^{(3)} &= w^{(4)}, \qquad & z^{(3)} &= Z^{(4)}, \quad & w^{(3)} &= Z^{(4)}W^{(4)}, \\
u^{(3)} &= u^{(4)}v^{(4)}, \quad & v^{(3)} &= v^{(4)}, \qquad & u^{(3)} &= U^{(4)}, \quad & v^{(3)} &= U^{(4)}V^{(4)}, 
\end{alignat*}
where the induced involution $\sigma$ maps $(z^{(4)}, w^{(4)}) \mapsto (-z^{(4)}-4 (w^{(4)})^{-2}, \, -w^{(4)})$ 
on $\C^2_{(z^{(4)},w^{(4)})}$ and $(u^{(4)}, v^{(4)}) \mapsto (-u^{(4)} + 4 (v^{(4)})^{-2}, \, -v^{(4)})$ on 
$\C ^2 _{(u^{(4)},v^{(4)})}$, respectively. 
In terms of the new charts there are two vertical leaves; the exceptional curve  
$\{w^{(4)} = 0\} \cup \{Z^{(4)} = 0\} \oplus \{v^{(4)} = 0\} \cup \{U^{(4)} = 0\}$ 
and the proper image $\{W^{(4)} = 0\} \oplus \{V^{(4)} = 0\}$ of $\{w^{(3)} = 0\} \oplus \{v^{(3)} = 0\}$, 
as well as a singular point $a^{(4)}_t = \{(z^{(4)}, w^{(4)}) = (0, 0) \} \oplus 
b^{(4)}_t = \{(u^{(4)}, v^{(4)}) = (0, 0) \}$. 
\par
{\bf 8.} Blowup at $a^{(4)}_t \oplus b^{(4)}_t$ produces new charts 
$(z^{(5)}, w^{(5)}) \oplus (u^{(5)}, v^{(5)})$ and 
$(Z^{(5)}, W^{(5)}) \oplus (U^{(5)}, V^{(5)})$ such that
\begin{alignat*}{4}
z^{(4)} &= z^{(5)}w^{(5)}, \quad & w^{(4)} &= w ^{(5)}, \qquad & z^{(4)} &=Z ^{(5)}, \quad & w^{(4)} &= Z^{(5)}W^{(5)}, \\
u^{(4)} &= u^{(5)}v^{(5)}, \quad & v^{(4)} &= v ^{(5)}, \qquad & u^{(4)} &=U ^{(5)}, \quad & v^{(4)} &= U^{(5)}V^{(5)},
\end{alignat*}
where the induced involution $\sigma$ maps $(z^{(5)}, w^{(5)}) \mapsto (z^{(5)}+ 4 (w^{(5)})^{-3}, \, -w^{(5)})$ 
on $\C ^2 _{(z^{(5)},w^{(5)})}$ and $(u^{(5)}, v^{(5)}) \mapsto (u^{(5)}- 4 (v^{(5)})^{-3}, \, -v^{(5)})$ on 
$\C ^2 _{(u^{(5)},v^{(5)})}$, respectively.   
In terms of the new charts there are two vertical leaves; the exceptional curve 
$\{w^{(5)} = 0\} \cup \{Z^{(5)} = 0\} \oplus \{v^{(5)} = 0\} \cup \{U^{(5)} = 0\}$ and 
the proper image $\{W^{(5)} = 0\} \oplus \{V^{(5)} = 0\}$ of $\{w^{(4)}=0\} \oplus \{v^{(4)} = 0\}$, 
as well as a singular point $a^{(5)}_t = \{(z^{(5)}, w^{(5)}) = (0, 0) \} \oplus  
b^{(5)}_t = \{(u^{(5)}, v^{(5)}) = (0, 0)\}$. 
\par
{\bf 9.} Blowup at $a^{(5)}_t \oplus b^{(5)}_t$ produces new charts 
$(z^{(6)}, w^{(6)}) \oplus (u^{(6)}, v^{(6)})$ and 
$(Z^{(6)}, W^{(6)}) \oplus (U^{(6)}, V^{(6)})$ such that 
\begin{alignat*}{4}
z^{(5)} &= z^{(6)}w^{(6)}, \quad & w^{(5)} &= w^{(6)}, \qquad & z^{(5)} &= Z^{(6)}, \quad & w^{(5)} &= Z^{(6)}W^{(6)}, \\
u^{(5)} &= u^{(6)}v^{(6)}, \quad & v^{(5)} &= v^{(6)}, \qquad & u^{(5)} &= U^{(6)}, \quad & v^{(5)} &= U^{(6)}V^{(6)},
\end{alignat*}
where the induced involution $\sigma$ maps $(z^{(6)}, w^{(6)}) \mapsto (-z^{(6)}-4(w^{(6)})^{-4}, \, -w^{(6)})$ 
on $\C ^2_{(z^{(6)},w^{(6)})}$ and $(u^{(6)}, v^{(6)}) \mapsto (-u^{(6)} +4 (v^{(6)})^{-4}, \, -v^{(6)})$ on 
$\C ^2 _{(u^{(6)},v^{(6)})}$, respectively. 
In terms of the new charts there are two vertical leaves; the exceptional curve  
$\{w^{(6)} = 0\} \cup \{Z^{(6)} = 0\} \oplus \{v^{(6)} = 0\} \cup \{U^{(6)} = 0\}$ and 
the proper image $\{W^{(6)} = 0\} \oplus \{V^{(6)} = 0\}$ of $\{w^{(5)} = 0\} \oplus \{v^{(5)} = 0\}$, 
as well as a singular point $a^{(6)}_t = \{(z^{(6)}, w^{(6)}) = (t/2, 0) \} \oplus 
b^{(6)}_t = \{(u^{(6)}, v^{(6)}) = (-t/2, 0)\}$. 
\par
{\bf 10.} Blowup at $a^{(6)}_t \oplus b^{(6)}_t$ produces new charts 
$(z^{(7)}, w^{(7)}) \oplus (u^{(7)}, v^{(7)})$ and 
$(Z^{(7)}, W^{(7)}) \oplus (U^{(7)}, V^{(7)})$ such that 
\begin{alignat*}{4}
z^{(6)} &= \phantom{-} t/2 + z^{(7)}w^{(7)}, \quad & w^{(6)} &= w^{(7)}, \qquad & 
z^{(6)} &= \phantom{-}t/2 +Z^{(7)}, \quad & w^{(6)} &= Z^{(7)}W^{(7)}, \\
u^{(6)} &=-t/2 + u^{(7)}v^{(7)}, \quad & v^{(6)} &= v^{(7)}, \qquad 
& u^{(6)} &=-t/2 +U^{(7)}, \quad & v^{(6)} &= U^{(7)}V^{(7)}, 
\end{alignat*}
where the induced involution $\sigma$ maps $(z^{(7)}, w^{(7)}) \mapsto 
(z^{(7)} + t/w^{(7)} + 4 (w^{(7)})^{-5}, \, -w^{(7)})$ on $\C ^2 _{(z^{(7)},w^{(7)})}$ and 
$(u^{(7)}, v^{(7)}) \mapsto (u^{(7)} -t/v^{(7)} - 4 (v^{(7)})^{-5}, \, -v^{(7)})$ on  
$\C ^2 _{(u^{(7)},v^{(7)})}$, respectively.    
In terms of the new charts there are two vertical leaves; the exceptional curve  
$\{w^{(7)} = 0\} \cup \{Z^{(7)} = 0\} \oplus \{ v^{(7)}=0\} \cup \{ U^{(7)}=0\}$ 
and the proper image $\{W^{(7)} = 0\} \oplus \{V^{(7)} = 0\}$ of $\{w^{(6)} = 0\} 
\oplus \{v^{(6)} = 0\}$, as well as a singular point  
$a^{(7)}_t = \{(z^{(7)}, w^{(7)}) = (1/2, 0)\} \oplus 
b^{(7)}_t = \{(u^{(7)}, v^{(7)}) = (1/2, 0)\}$. 
\par
{\bf 11.} Blowup at $a^{(7)}_t \oplus b^{(7)}_t$ produces new charts 
$(z^{(8)}, w^{(8)}) \oplus (u^{(8)}, v^{(8)})$ and 
$(Z^{(8)}, W^{(8)}) \oplus (U^{(8)}, V^{(8)})$ such that
\begin{alignat*}{4}
z^{(7)} &= 1/2 + z^{(8)}w^{(8)}, \quad & w^{(7)} &= w^{(8)}, \quad & z^{(7)} &= 1/2 + Z^{(8)}, \quad & w^{(7)} &= Z^{(8)}W^{(8)}, \\
u^{(7)} &= 1/2 + u^{(8)}v^{(8)}, \quad & v^{(7)} &= v^{(8)}, \quad & u^{(7)} &= 1/2 + U^{(8)}, \quad & v^{(7)} &= U^{(8)}V^{(8)}, 
\end{alignat*}
where the induced involution $\sigma$ maps $(z^{(8)}, w^{(8)}) \mapsto 
(-z^{(8)}-t (w^{(8)})^{-2} - 4 (w^{(8)})^{-6}, \, -w^{(8)})$ on $\C ^2 _{(z^{(8)},w^{8)})}$ and 
($u^{(8)}, v^{(8)}) \mapsto (-u^{(8)} + t (v^{(8)})^{-2} + 4 (v^{(8)})^{-6}, \, -v^{(8)})$ on 
$\C ^2 _{(u^{(8)},v^{(8)})}$, respectively.  
In terms of the new charts there is only one vertical leaf; the proper image  
$\{W^{(8)} = 0\} \oplus \{V^{(8)} = 0\}$ of $\{w^{(7)} = 0\} \oplus \{v^{(7)} = 0\}$.   
Observe that the exceptional curve $\{Z^{(8)} = 0\} \cup \{w^{(8)} = 0\} \oplus 
\{U^{(8)} = 0\} \cup \{v^{(8)} = 0\}$ is {\sl not} a vertical leaf and there 
is {\sl no} singular point of the foliation. 
\par
{\bf 12.} Composition of steps 6--11 leads to a proper modification 
$(W,p,p_{\pm}) \leftarrow (X,p,E_{\pm})$. 
The rest is just as mentioned in \S\ref{sec:intro}. 
Make a gluing $F = (S \setminus C) \cup (X/\sigma)$ to have a compact 
surface $F$ with an $A_1$-singularity $p \in F$; take its minimal resolution 
to get a smooth compact space $\overline{E}_t$; and finally remove the 
vertical leaves $V_t$ to obtain $E_t$. 
All these procedures are symbolically represented by Fig. \ref{fig:blowup7}. 
In order to make the final result exactly symplectic, we define the final 
charts $(z,w)$ and $(u,v)$ by $z^{(8)} = - z/2$, $w^{(8)} = w$; 
$u^{(8)} = -u/2$, $v^{(8)} = v$. 
Among all charts of $\overline{E}_t$ we have constructed, 
those which are disjoint with $V_t$ are exactly $\C^2_{(x,y)}$, 
$\C^2_{(z,w)}$ and $\C^2_{(u,v)}$. 
These three make an orbifold symplectic atlas of $E_t$. 
A careful check of steps 1--11 yields the desired relations 
\eqref{eqn:(z,w)-(u,v)}, \eqref{eqn:x-w}, \eqref{eqn:x-v} as well 
as formula \eqref{eqn:sigma2} for the involution $\sigma$.  
\par
It might be fun to think of Fig. \ref{fig:blowup7} as a ``deer" whose 
dual ``horns" $E_{\pm}$ are identified by $\sigma$, and whose ``nose" 
is just the fixed point $p \in X$ or the $A_1$-singularity 
$p \in X/\sigma$ arising from it.  
\section{Holomorphic Functions} \label{sec:hol}
We prove Theorem \ref{thm3}. 
Fixing $t \in T$ we do not refer to the dependence upon $t$.  
Any function holomorphic on $E_t$ and meromorphic on $\overline{E}_t$ 
is represented by a triple $\{H, K, L \}$ of functions $H = H(x,y)$, 
$K = K(z,w)$ and $L = L(u,v)$ entire in their respective variables 
such that $H = K = L$ under transformations \eqref{eqn:(z,w)-(u,v)}, 
\eqref{eqn:x-w} and \eqref{eqn:x-v}, as well as 
$K \circ \sigma = K$ and $L \circ \sigma = L$.  
\begin{lemma} \label{lem:polyn} 
We have $H \in \C[x, y]$, $K\in \C[z,w]$ and $L \in \C[u,v]$.  
\end{lemma}
{\it Proof}. As an entire holomorphic function of $(x, y)$, 
$H$ admits a Taylor expansion 
\[
H = \sum_{i,j=0}^{\infty} c_{ij} \, x^i y^j, \qquad c_{ij} \in \C. 
\]
Using relation \eqref{eqn:hirzebruch} we can rewrite it in terms of 
$(q_2, p_2)$ to have
\[
H = \sum_{i,j=0}^{\infty} (-1)^j c_{ij} \, q_2^{-(i+2j)} p_2^{-j}. 
\]
In order for this to be meromorphic on the vertical leaf $\{q_2 = 0\}$, 
there must be a nonnegative integer $N$ such that $c_{ij} = 0$ for every 
$i + 2 j > N$, which forces $H \in \C[x, y]$. 
Next we show that $K \in \C[z, w]$.  
Under transformation \eqref{eqn:x-w} we have $K = H \in \C[x, y]$ so 
that $K \in \C[z, w, w^{-1}]$. 
On the other hand, $K \in \C\{z, w\}$ (the convergent power series ring), 
since it is an entire function of $(z, w)$.  
Thus we have $K \in \C[z, w, w^{-1}] \cap \C\{z, w\} = \C[z, w]$. 
Similarly, $L \in \C[u, v]$. \hfill $\Box$
\par\medskip
Let us discuss the problem of finding a function $K$ such that  
\begin{equation} \label{eqn:inv-inv}
K \circ \sigma = K, \qquad K = K(z,w) \in \C[z,w],  
\end{equation}
where $\sigma : \C^2_{(z,w)} \circlearrowleft$ is the involution 
defined by formula \eqref{eqn:sigma2}. 
We begin by a simple reduction.  
Consider the decomposition of $K$ into even and odd components with 
respect to $w$: 
\begin{equation} \label{eqn:eod}
K = K^+ + K^-, \qquad K^{\pm} := (K \pm \check{K})/2,  
\qquad \check{K}(z,w) := K(z,-w). 
\end{equation}
\begin{lemma} \label{lem1}
If $K$ is a solution to problem \eqref{eqn:inv-inv}, then so are 
$K^{\pm}$.  
\end{lemma}
{\it Proof}. 
It suffices to show that if $K$ is a solution to problem 
\eqref{eqn:inv-inv} then so is $\check{K}$.  
First, it is obvious that $K \in \C[z,w]$ implies 
$\check{K} \in \C[z,w]$. 
Next, observe that 
\[
\begin{split}
(\check{K}\circ\sigma)(z,w) &= \check{K}(8 w^{-6} + 2 t w^{-2} - z, \, -w) 
= K(8 w^{-6} + 2 t w^{-2} - z, \, w) \\
&=  K(8 (-w)^{-6} + 2 t (-w)^{-2} - z, \, -(-w)) = (K \circ \sigma)(z,-w) \\ 
&= K(z,-w) = \check{K}(z,w),  
\end{split}
\]
where $K \circ \sigma = K$ is used in the fifth equality. 
Thus $\check{K}$ is also a solution to problem \eqref{eqn:inv-inv}. 
\hfill $\Box$ \par\medskip 
Let $K$ be a solution to problem \eqref{eqn:inv-inv} 
and put $\xi = w^2$.  
The even component of $K$ can be written $K^+(z,w) = F(z, w^2)$ 
with $F = F(z, \xi)$ being a solution to the problem   
\begin{equation} \label{eqn:e-inv-inv}
F \circ \tau = F, \qquad F = F(z,\xi) \in \C[z,\xi],  
\end{equation}
where $\tau : \C^2_{(z,\xi)} \circlearrowleft$ is an involution 
$(z,\xi) \mapsto \left(8 \xi^{-3} + 2 t \xi^{-1} - z, \, \xi \right)$.   
There is a particular solution $E(z,\xi) := z (\xi^3 z - 2t \xi^2 - 8)$ 
to problem \eqref{eqn:e-inv-inv}, which plays an important role 
in the following.  
\begin{lemma} \label{lem2}  
Any nontrivial solution to problem \eqref{eqn:e-inv-inv} must 
be of the form 
\begin{equation} \label{eqn:e-sol}
F(z,\xi) = \sum_{m=0}^M f_m(\xi) E^m(z, \xi), \qquad  
f_m(\xi) \in \C[\xi],  
\end{equation}
where $M \ge 0$ and $f_M(\xi)$ is a nonzero polynomial of $\xi$.  
\end{lemma}
{\it Proof}. 
We prove the lemma by induction on $\deg_z F(z,\xi)$.  
If $\deg_z F(z,\xi) = 0$ then formula \eqref{eqn:e-sol} obviously holds 
with $M = 0$.  
Suppose that $\deg_z F(z,\xi) \ge 1$ and put $f_0(\xi) := F(0,\xi) \in \C[\xi]$. 
Notice that $F_0(z,\xi) := F(z,\xi) - f_0(\xi) \in \C[z,\xi]$ is also a 
solution to problem \eqref{eqn:e-inv-inv}. 
It is divisible by $z$, that is, $F_0(z,\xi) = z F_1(z,\xi)$ for some 
$F_1(z,\xi) \in \C[z,\xi]$. 
The $\tau$-invariance of $F_0(z, \xi)$ implies 
$z F_1(z, \xi) = (8 \xi^{-3} + 2 t \xi^{-1} - z) F_2(z, \xi)$, where 
$F_2(z, \xi) := F_1(8 \xi^{-3} + 2 t \xi^{-1} - z, \xi) \in \C(\xi)[z]$. 
Writing $F_2(z,\xi) = - F_3(z, \xi)/a(\xi)$ with $F_3(z, \xi) \in \C[z, \xi]$ 
and $a(\xi) \in \C[\xi]$, we obtain $a(\xi) \xi^3 z F_1(z,\xi) = 
(\xi^3 z-8-2 t \xi^2) F_3(z, \xi)$ in $\C[z, \xi]$. 
Since the right-hand side is divisible by $\xi^3 z-8-2 t \xi^2$, 
so must be the left-hand side, but $a(\xi) \xi^3 z$ and 
$\xi^3 z-8-2 t \xi^2$ have no common factor, so that $F_1(z, \xi)$ 
must be divisible by $\xi^3 z-8-2 t \xi^2$, that is, 
$F_1(z, \xi) = (\xi^3 z-8-2 t \xi^2) F_4(z, \xi)$ for 
some $F_4(z,\xi) \in \C[z, \xi]$. 
Thus we have $F_0(z, \xi) = E(z,\xi) F_4(z, \xi)$. 
Since $F_0(z,\xi)$ and $E(z,\xi)$ are $\tau$-invariant, 
$F_4(z, \xi)$ is also $\tau$-invariant and hence yields a solution 
to problem \eqref{eqn:e-inv-inv} with $\deg_z F_4(z, \xi) = 
\deg_z F_0(z, \xi) -2 = \deg_z F(z, \xi) -2$. 
By induction hypothesis we can write $F_4(z,\xi) = \sum_{m=1}^M 
f_m(\xi) E^{m-1}(z, \xi)$ for some $f_m(\xi) \in \C[\xi]$.   
Substituting this into $F(z,\xi) = f_0(\xi) + E(z, \xi) 
F_4(z, \xi)$ yields formula \eqref{eqn:e-sol}. 
The induction is complete.  \hfill $\Box$ \par\medskip 
On the other hand, the odd component of $K$ can be written $K^-(z,w) = 
w G(z,w^2)$ with $G = G(z, \xi)$ being a solution to the problem  
\begin{equation} \label{eqn:o-inv-inv}
G \circ \tau = - G, \qquad G = G(z,\xi) \in \C[z,\xi].   
\end{equation}
Notice that $\varDelta(z,\xi) := \xi^3 z - t \xi^2 - 4$ is a particular 
solution to problem \eqref{eqn:o-inv-inv}.  
\begin{lemma} \label{lem3} 
Any nontrivial solution to problem \eqref{eqn:o-inv-inv} must 
be of the form 
\begin{equation} \label{eqn:o-sol}
G(z,\xi) = \varDelta(z, \xi) \sum_{n=0}^N g_n(\xi) E^n(z, \xi), \qquad  
g_n(\xi) \in \C[\xi],   
\end{equation}
where $N \ge 0$ and $g_N(\xi)$ is a nonzero polynomial of $\xi$.  
\end{lemma}
{\it Proof}.  
Substituting $z = t \xi^{-1} + 4 \xi^{-3}$ into the skew 
$\tau$-invariance $G(8 \xi^{-3} + 2 t \xi^{-1} - z, \xi) = - G(z,\xi)$ 
yields $G(t \xi^{-1} + 4 \xi^{-3}, \xi) = - 
G(t \xi^{-1} + 4 \xi^{-3}, \xi)$, which forces 
$G(t \xi^{-1} + 4 \xi^{-3}, \xi) = 0$. 
Thus $G(z, \xi)$ is divisible by $z - t \xi^{-1} - 4 \xi^{-3}$ in 
$\C(\xi)[z]$, that is, $G(z, \xi) = (z - t \xi^{-1} - 4 \xi^{-3}) 
G_1(z,\xi)$ for some $G_1(z, \xi) \in \C(\xi)[z]$. 
Writing $G_1(z, \xi) = G_2(z, \xi)/b(\xi)$ with $G_2(z, \xi) \in 
\C[z, \xi]$ and $b(\xi) \in \C[\xi]$, we obtain 
$b(\xi) \xi^3 G(z, \xi) = \varDelta(z,\xi) G_2(z, \xi)$ in $\C[z, \xi]$. 
Since the right-hand side is divisible by $\varDelta(z,\xi)$, 
so must be the left-hand side, but $b(\xi) \xi^3$ and 
$\varDelta(z,\xi)$ have no common factor, so that $G(z, \xi)$ 
must be divisible by $\varDelta(z,\xi)$, that is, 
$G(z, \xi) = \varDelta(z,\xi) G_3(z, \xi)$ for some 
$G_3(z,\xi) \in \C[z, \xi]$. 
Since $G(z,\xi)$ and $\varDelta(z,\xi)$ are skew $\tau$-invariant, 
$G_3(z, \xi)$ is $\tau$-invariant and so yields a solution 
to problem \eqref{eqn:e-inv-inv}. 
Lemma \ref{lem2} then allows us to write 
$G_3(z,\xi) = \sum_{n=0}^N g_n(\xi) E^n(z,\xi)$ for some $g_n(\xi) 
\in \C[\xi]$, which leads to representation \eqref{eqn:o-sol}.   
\hfill $\Box$
\par\medskip
Now a general solution to problem \eqref{eqn:inv-inv} can be written  
$K(z,w) = K^+(z,w)+ K^-(z,w)$ with $K^+(z,w) = F(z, w^2)$, 
$K^-(z,w) = w G(z, w^2)$, where $F(z,\xi)$ and 
$G(z,\xi)$ are as in formulas \eqref{eqn:e-sol} and 
\eqref{eqn:o-sol} respectively. 
Recall that we have $H(x,y) = K(z,w)$ under transformation 
\eqref{eqn:x-w}.  
Let $H(x,y) = H^+(x,y) + H^-(x,y)$ be the decomposition 
parallel to the one $K(z,w) = K^+(z,w)+ K^-(z,w)$. 
Notice that $H^{\pm}(x,y) \in \C[x,x^{-1},y]$. 
Observe that  
\begin{equation} \label{eqn:ED}
E(z, w^2) = 4 y^2 + 4 x^{-1}y + x^{-2} - x^{-1}(4x^2+t)^2, \qquad 
w \varDelta(z, w^2) = 2x^{-2} y + x^{-3},   
\end{equation}
under relation \eqref{eqn:x-w}. 
Indeed, the second formula readily follows from \eqref{eqn:x-w}, 
while the first formula is derived from the second one and the relation   
$E(t,w^2) = w^{-8} \{w \varDelta(z, w^2) \}^2 - w^2 \left(t + 4 w^{-4}\right)^2$.  
Formulas \eqref{eqn:ED} are substituted into formulas \eqref{eqn:e-sol} 
and \eqref{eqn:o-sol} to find 
\begin{equation} \label{eqn:Hpm}
\begin{split}
H^+(x,y) &= \sum_{m=0}^M f_m(x^{-1})
 \left\{4 y^2 + 4 x^{-1}y + x^{-2} - x^{-1}(4x^2+t)^2\right\}^m, \\
H^-(x,y) &= \left(2x^{-2} y + x^{-3}\right) 
\sum_{n=0}^N g_n(x^{-1}) 
\left\{4 y^2 + 4 x^{-1}y + x^{-2} - x^{-1}(4x^2+t)^2\right\}^n,   
\end{split}
\end{equation}
with $f_m(\xi) \in \C[\xi]$ and $g_n(\xi) \in \C[\xi]$, where if 
$H^+(x,y)$ is nontrivial then $M \ge 0$ and $f_M(\xi)$ is a nonzero 
polynomial, while if $H^-(x,y)$ is nontrivial then $N \ge 0$ and 
$g_N(\xi)$ is a nonzero polynomial. 
By convention we put $M = -1$ resp. $N = -1$ if $H^+(x,y)$ resp. 
$H^-(x,y)$ is trivial. 
Suppose that $H(x,y)$ is nontrivial, so that at least one of 
$H^{\pm}(x,y)$ is nontrivial. 
\begin{lemma} \label{lem4}
We have $M > N$ and $H^+(x,y)$ must be nontrivial. 
\end{lemma}
{\it Proof}.  
Expanding formulas \eqref{eqn:Hpm} into powers of $y$ yield  
\begin{equation} \label{eqn:Hpm-y} 
\begin{split}
H^+(x,y) &= 2^{2M} f_M(x^{-1}) \, y^{2M} + 
2^{2M} M x^{-1} \, f_M(x^{-1}) \, y^{2M-1} + \cdots, \\
H^-(x,y) &= 2^{2N+1} x^{-2} g_N(x^{-1}) \, y^{2N+1} + \cdots, 
\end{split}
\end{equation}
where $\cdots$ denotes lower-degree terms with respect to $y$.  
If $M \le N$ then $2M < 2N+1$ in formula \eqref{eqn:Hpm-y} so that 
$H(x,y) =  2^{2N+1} x^{-2} g_N(x^{-1}) y^{2N+1} + \cdots \in \C[x,y]$ 
and hence $x^{-2} g_N(x^{-1}) \in \C[x]$. 
This is possible only if $g_N(\xi)$ is the zero polynomial, in which 
case $H^-(x,y)$ is trivial with $N = -1$; then $M = -1$ and so 
$H^+(x,y)$ is also trivial. 
This contradiction shows that $M > N \ge -1$. 
Since $M$ is nonnegative, $H^+(x,y)$ must be nontrivial. \hfill $\Box$
\begin{lemma} \label{lem5} 
We have $M = 0$, $N = -1$ and $H(x,y) = c \in \C^{\times}$.  
\end{lemma}
{\it Proof}. By Lemma \ref{lem4} we have $M > N$ and hence $2M > 2N+1$, 
so that formulas \eqref{eqn:Hpm-y} yield $H(x,y) = 2^{2M} f_M(x^{-1}) \, 
y^{2M} + \cdots \in \C[x,y]$, which implies $f_M(x^{-1}) \in \C[x]$. 
On the other hand we have $f_M(x^{-1}) \in \C[x^{-1}]$. 
Thus $f_M(\xi)$ must be a constant, say, $c \in \C$. 
Since $H^+(x,y)$ is nontrivial, $f_M(\xi) = c \in \C^{\times}$ must 
be a nonzero constant.  
To show that $M = 0$, suppose the contrary that $M \ge 1$. 
In the first case where $N = M-1$, formulas \eqref{eqn:Hpm-y} imply 
\[
H(x,y) = 2^{2M} c \, y^{2M} + \{(M \cdot 2^{2M} \cdot c) x^{-1} 
+ 2^{2M-1} x^{-2} g_{M-1}(x^{-1})\} y^{2M-1} + \cdots \in \C[x,y],  
\]
and so $(M \cdot 2^{2M} \cdot c) x^{-1} + 2^{2M-1} x^{-2} 
g_{M-1}(x^{-1}) \in \C[x]$, which is impossible. 
In the second case where $N < M-1$, formulas \eqref{eqn:Hpm-y} imply 
$H(x,y) = 2^{2M} c \, y^{2M} + (M \cdot 2^{2M} \cdot c) x^{-1} y^{2M-1} + 
\cdots \in \C[x,y]$, and so $(M \cdot 2^{2M} \cdot c) x^{-1} \in \C[x]$, 
which is also impossible. 
Thus $M = 0$ and $H^+(x,y) = c$. 
Since $N < M = 0$, we have $N = -1$ and $H^-(x,y) = 0$ so that $H(x,y) = c$.  
\hfill $\Box$ \par\medskip
With the proof of Lemma \ref{lem5} above, Theorem \ref{thm3} has also been 
established completely. 
\par
Formula \eqref{eqn:Hpm} can be used to construct a meromorphic 
function on $E_t$ that is holomorphic on $E_t \setminus \{x=0\}$ with 
poles only along $\{x = 0\}$. 
There is a connection of this formula with a proof of the Painlev\'e 
property.   
In a qualitative proof of it, which does not use isomonodromic deformations, 
it is crucial to deal with a kind of Lyapunov function that can control 
the trajectories near the vertical leaves. 
As a Lyapunov function for $\rP_{\mathrm{I}}$ we usually employ    
\[
U(x,y,t) = 2 H_{\mathrm{I}}(x,y,t) + \dfrac{y}{x} = y^2 - 4 x^3 - 2t x + 
\dfrac{y}{x} 
\]
as in \cite[formula (5)]{OT} or \cite[formula (3.8)]{Shimomura}.   
This function is just obtained by putting $M = 1$, $f_1(\xi) = 1/4$, 
$f_0(\xi) = \xi(t^2-\xi)/4$ and $N = -1$, i.e., $H^-(x,y,t) = 0$ 
in formula \eqref{eqn:Hpm}. 
\par
A quite different proof, but still of a qualitative nature, for the 
Painlev\'e property has been proposed by H.~Chiba in his framework of 
Painlev\'e equations on weighted projective spaces \cite{Chiba}.    

\end{document}